\documentclass[reqno,12pt]{amsart}
\usepackage{amssymb,graphicx}
\def\N{{\mathbb N} } 
\def\R{{\mathbb R} } 
\def\jj{{\bf j} }    \def\ii{{\bf i} }
\def\JJ{{\mathcal J} }  \def\HH{{\mathcal H} }
\def\diam{{\rm diam\,}}  \def\dist{{\rm dist\,}} \def\inte{{\rm int\,}}
\def\conv{{\rm conv\,}}  \def\supp{{\rm supp\,}}

\newtheorem{Thm}{Theorem} %[section]

\newtheorem{Prop}[Thm]{Proposition}
\newtheorem{lem}[Thm]{Lemma}

\newtheorem{Exa}[Thm]{Example}
\newtheorem{rem}[Thm]{Remark}

\newenvironment{Proof}[0]
{\medskip \noindent {\it Proof.} \ }{\ \hfill $\Box$\medskip}

\begin{document}
\title  {Differentiability of fractal curves}
\author{Christoph Bandt and Alexey Kravchenko}
\thanks{A.K. was supported by grant RFBR 09-01-90202-Mong\_a and grant ADTP
``Development of Scientific Potential of Higher Education''
2.1.1/3707 (Russia). Most work was done when A.K. was employed at
Greifswald University.}

\begin{abstract} While self-similar sets have no tangents
at any single point, self-affine curves can be smooth. We consider
plane self-affine curves without double points and with two
pieces. There is an open subset of parameter space for which the
curve is differentiable at all points except for a countable set.
For a parameter set of codimension one, the curve is continuously
differentiable. However, there are no twice differentiable
self-affine curves in the plane, except for parabolic arcs.
\end{abstract}
\maketitle
\noindent MSC classification: Primary 28A80, Secondary 26A27, 28A75 \vspace{12mm}\\
Christoph Bandt, Alexey Kravchenko\\
Institute for Mathematics and Informatics\\
Arndt University\\ 17487 Greifswald, Germany\\
e-mail: {bandt@uni-greifswald.de}, {alexey14@gmail.com}
\pagebreak

%%%%%%%%%%%%%%%%%%%%
\section{Overview}
{\bf I. } It is well-known that fractals are not differentiable.
Nevertheless some fractals, like the Appolonian gasket, or
parabolic Julia sets, possess tangents at some of their points.
For self-similar sets
$$A=f_1(A)\cup ...\cup f_m(A)  $$
where the $f_i$ are contracting similarity maps on $\R^n$
\cite{Hut,Fal,Bar}, there are no such exception points:

\begin{Thm}
(i) If $A$ is a self-similar set which spans $\R^n$ and $x\in A ,$
there does not exist a tangent hyperplane of $A$ at $x.$\\
(ii) If $\mu$ is a self-similar measure and $x\in \supp\mu,$ there
does not exist an approximate tangent hyperplane of $\mu$ at $x$
in the measure-theoretical sense.
\end{Thm}

Precise definitions and the proof are given in Section 2. Related
results on self-conformal sets with separation condition can be
found in K\"aenm\"aki \cite{K0,K1} and the references given there.
Theorem 1 directly extends to self-conformal sets without
separation condition - see Remark 6.\vspace{2.5ex}

{\bf II. } Self-affine curves, however, can be smooth, as shown in
\cite{K}. Here we prove that differentiability of self-affine
curves is not an exception, but a rather generic phenomenon. We
shall consider self-affine curves $J$ in the plane with two
pieces:
$$J=f_1(J)\cup f_2(J)  $$
where $f_1,f_2$ are contracting affine maps in $\R^2$ with
positive eigenvalues and with fixed points $e_1,e_2$ respectively,
and $f_2(e_1)=f_1(e_2).$ In Section 3 we show that the structure
of $J$ is determined, up to an affine coordinate transformation,
by the eigenvalues $\lambda_1,\nu_1$ of $f_1$ and
$\lambda_2,\nu_2$ of $f_2.$ We take $\lambda_i\ge\nu_i>0,$ exclude
similarity maps, and assume that the eigenvectors associated with
$\lambda_1$ and $\lambda_2$ do not coincide. In Section 3, these
conditions are stated more technically as (1) and (2).

\begin{Thm}
(i) Under the above assumptions, the curve $J$ is differentiable
at all points $x\in J$ except for a countable set if \[
\lambda_1+\nu_2<1\quad \mbox{ and }\quad \lambda_2+\nu_1<1\  . \]
(ii) If this condition holds, the curve is continuously
differentiable if and only if the one-sided tangents at the
intersection point $z=f_2(e_1)$ coincide:
\[ \nu_1\nu_2= (1-\lambda_1-\nu_2)(1-\lambda_2-\nu_1)\,\ . \]
\end{Thm}

Figure 1 shows examples of everywhere and almost everywhere
differentiable self-affine curves. They were standardized to
represent functions
\[ x_2=\psi(x_1)\quad\mbox{ with }\quad
\psi(-1)=\psi(1)=1\quad\mbox{ and }\quad \psi'(-1)=-1,\
\psi'(1)=1\, .
\]
According to Theorem 2, there is a four-parameter family of such
functions which are almost everywhere differentiable, and there is
a three-parameter subfamily of continuously differentiable
self-affine functions. In Section 3 we discuss the properties of
such curves. In Section 4 we prove Theorem 2.\vspace{7.5ex}

\begin{figure}
\includegraphics[width=.32\textwidth]{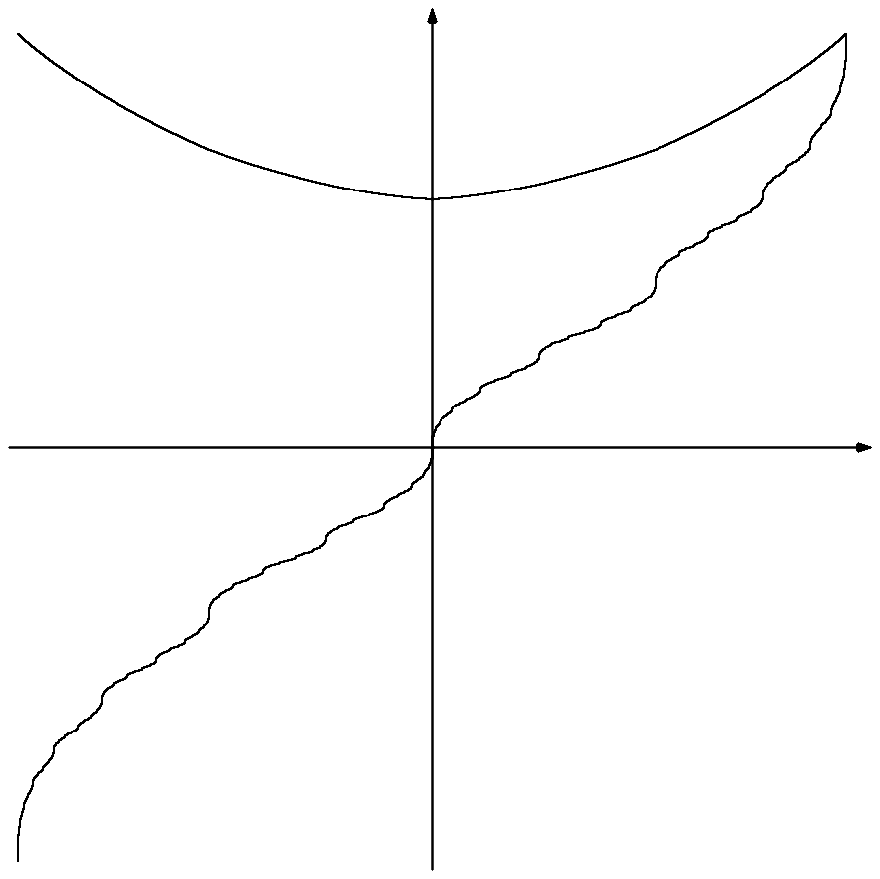}
\includegraphics[width=.32\textwidth]{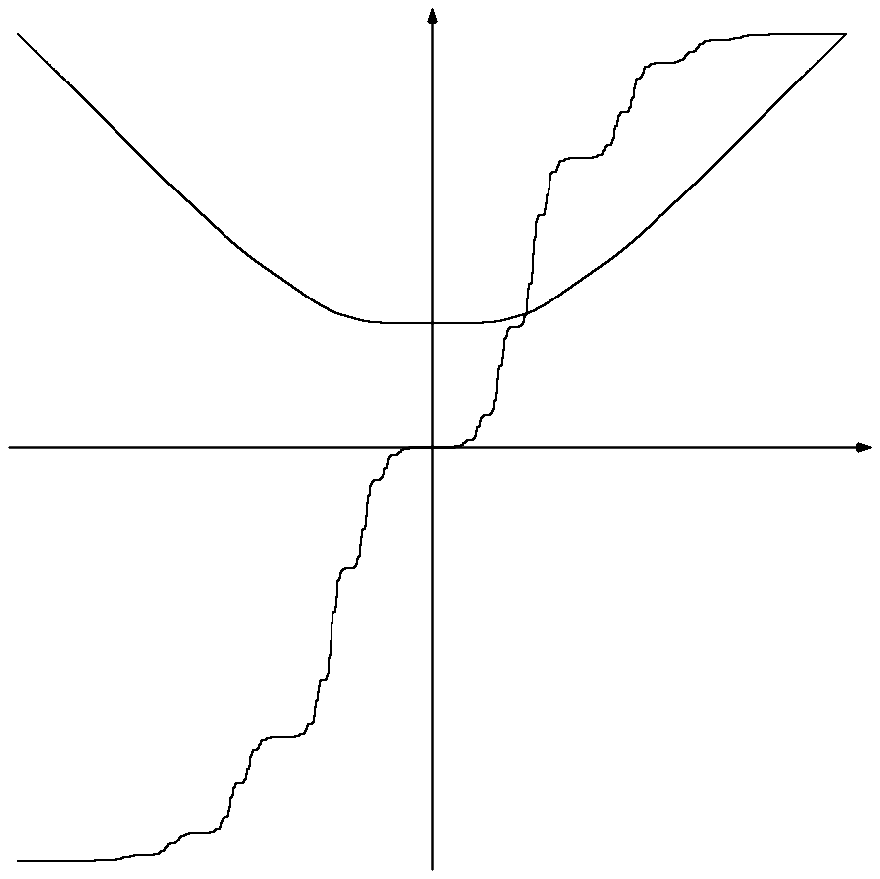}
\includegraphics[width=.32\textwidth]{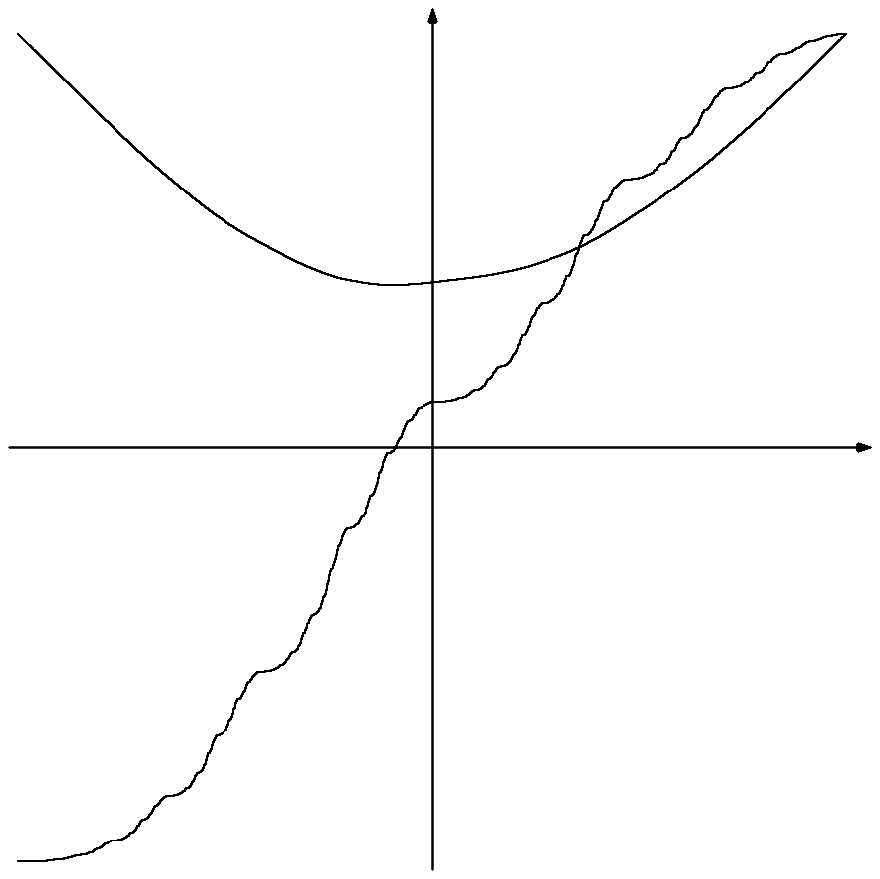}\vspace{3.5ex}\\
\includegraphics[width=.32\textwidth]{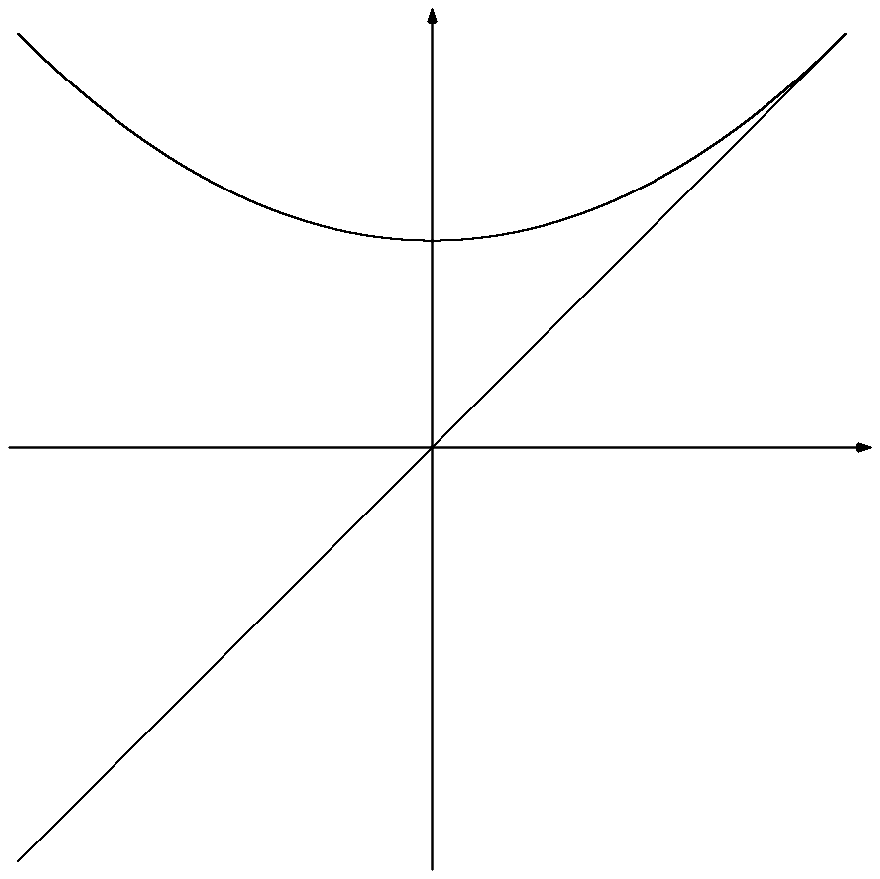}
\includegraphics[width=.32\textwidth]{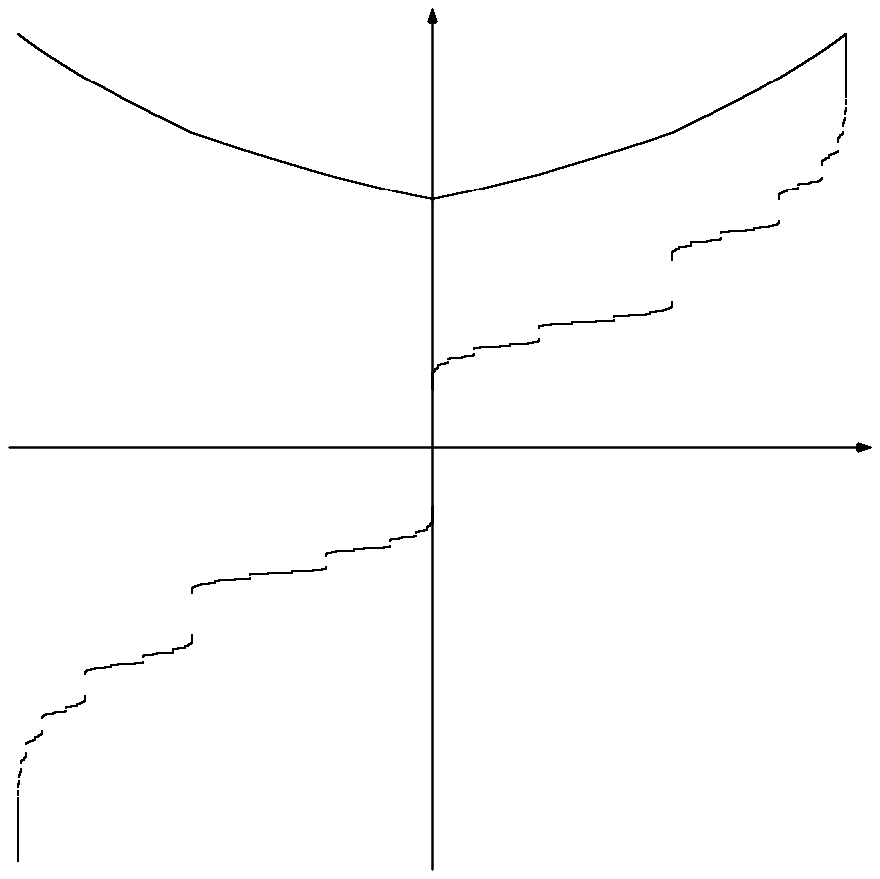}
\includegraphics[width=.32\textwidth]{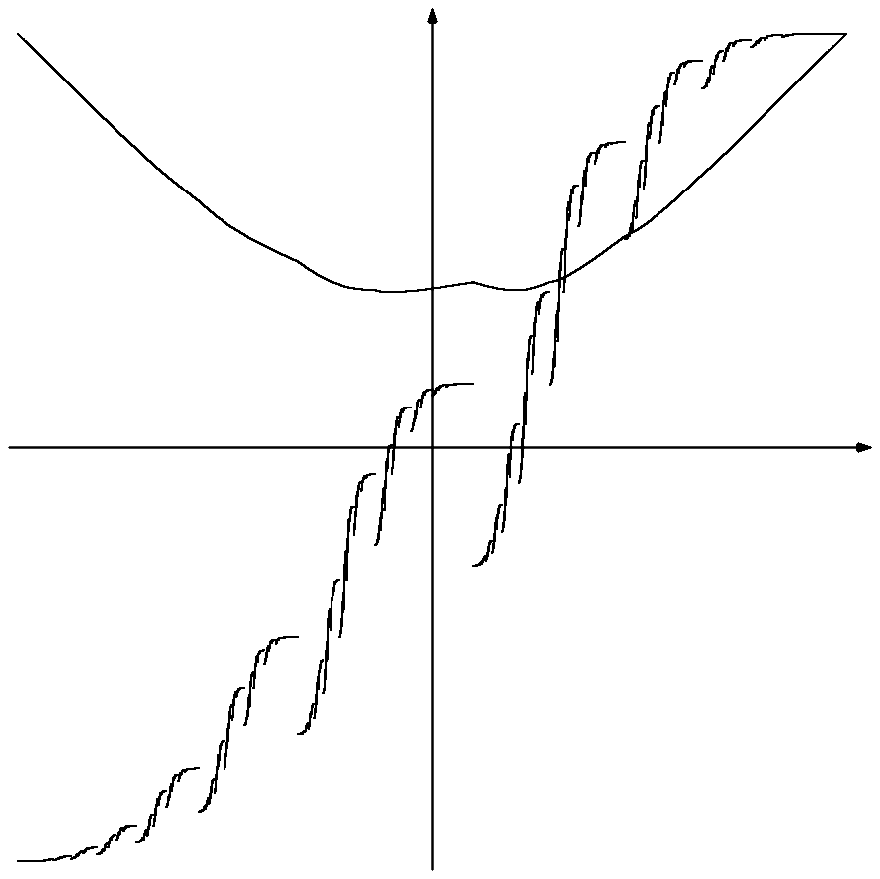}\vspace{1.5ex}\\
\caption{Upper row: several smooth self-affine functions with
their derivatives. Lower row: parabola and two self-affine
functions which are differentiable up to a countable set of points
where only one-sided derivatives exist}
\end{figure}

{\bf III. } Nevertheless, self-affine curves are not very smooth.
With exception of segment and parabola, they are not $C^2.$ The
following theorem is proved in Section 5.

\begin{Thm}
(i) If a plane curve $J,$ parametrized as $x=\phi(t),
\phi:[0,1]\to J,$ is two times continuously differentiable with
$\phi(0)=0$ and $\phi''(0)\not= 0,$ and there exists a contractive
linear map $h\not= 0$ which maps $J$ into itself, then $J$ is a
parabolic arc.\\
(ii) Except for parabolic arcs and segments, there are no twice
continuously differentiable self-affine curves in the plane.
\end{Thm}

\section{No tangents to self-similar sets}
An affine subspace $V$ of $\R^n$ is {\it tangent} to a set $B$ at
a point $x\in V\cap B$ if for each $\delta>0$  there is an
$\epsilon>0$ such that \[\frac{\dist (b,V)}{|x-b|}<\delta \mbox{
for all } b\in B\setminus\{ x\} \mbox{ with }|x-b|<\epsilon\ .\]
Here $\dist(b,V)=\inf \{ |b-v|\,| v\in V\}.$ In other words, all
chords from $x$ to a point $b\in B\cap U_\epsilon (x)$ should
subtend an angle smaller $\alpha=\arcsin\delta$ with $V.$

\begin{lem} Let $A$ be a bounded set which spans $\R^n,$ and let $a_0,a_1,...,a_n$ be $n+1$
points of $A$ in general position. Then there is an $\phi>0$ such
that for all $y\in A$ and each hyperplane $V$ through $y,$ there
exists $i\in\{ 0,1,...,n\}$ such that the chord $ya_i$ forms an
angle $\ge\phi$ with $V.$ \end{lem}

{\it Proof. } Let $B$ be a closed ball with radius $R$ which lies
within the simplex $S=\conv\{ a_0,..,a_n\}.$ Now take an arbitrary
point $y\in A$ and an arbitrary hyperplane $V$ through $y.$ Let
$W=\{ x\in\R^n\,|\, \dist (x,V)<R\} .$ Then $B$ is not a subset of
$W.$ Since $W$ is convex, and $B$ lies in the convex hull of the
$a_j,$ it follows that there is at least one point $a_i$ outside
$W.$ For the angle $\beta$ between $ya_i$ and $V$ we have
\[ \sin \beta = \frac{\dist (a_i,V)}{|a_i-y|} \ge \frac{R}{\diam A}\  .\]
We proved the lemma for $\phi=\arcsin \frac{R}{\diam A}.$ \hfill
$\Box$\medskip

{\it Proof of Theorem 1 (i). } The lemma prevents approximation by
a plane on global scale. We show that for a {\it self-similar} set
$A,$ Lemma 4 remains true in small neighborhoods. Let $\delta=\sin
\frac{\phi}{2},$ and let $x\in A$ and $\epsilon>0$ be arbitrarily
chosen. There is a small piece $A_\jj=f_\jj(A)$ of $A$ within
$U_\epsilon(x).$ Here $\jj=j_1...j_k$ and
$f_\jj=f_{j_1}...f_{j_k}.$ The similarity map $f_\jj^{-1}$ maps
$A_\jj$ onto $A,$ and any hyperplane $V'$ through $x$ onto a
corresponding hyperplane $V$ through $y=f_\jj^{-1}(x).$ Lemma 4
says that there is a chord $ya_i$ which subtends an angle
$\ge\phi$ with $V.$ Since the chord $xf_\jj(a_i)$ within
$U_\epsilon(x)$ subtends the same angle with $V',$ the definition
of tangent is not fulfilled for $V'$ and $x.$ So there is no
tangent hyperplane to $A$ in $x.$\hfill $\Box$\medskip

To deal with tangents of measures, we have to specify concepts.
For a finite Borel measure $\mu,$ the {\it support} $\supp\mu$
contains all those points $y$ for which $\mu(U_\epsilon(y))>0$ for
all $\epsilon>0.$ An affine subspace $V$ of $\R^n$ is an {\it
approximate tangent} to $\mu$ at a point $x\in V\cap\supp\mu$ if
for each $\delta>0$ there is an $\epsilon>0$ with
\[ \mu\{ b\in U_\epsilon(x)\,|\,
\frac{\dist(b,V)}{|x-b|}>\delta\}\, / \,\mu (U_\epsilon(x)) <
\delta\ .\] Now we allow for chords $xb$ which subtend an angle
$>\alpha =\arctan\delta$ with $V.$ But the percentage of endpoints
$b\in U_\epsilon(x)$ with this property, measured by $\mu ,$
should converge to zero with $\epsilon .$ See Mattila \cite{Mat},
Chapter 15, for related concepts.\vspace{1ex}

We are going to derive Theorem 1 (ii) from the above proof of (i).
First we reformulate a global fact for arbitrary measures, giving
a lower bound $\eta$ for the number of exceptions.

\begin{lem} Let $\mu$ be a probability measure on $\R^n$ such that
$A=\supp\mu$ is bounded and contains $n+1$ points $a_0,...,a_n$ in
general position. Then there are $\psi>0, \gamma>0, \eta>0$ such
that for all $y\in \supp\mu$ and each hyperplane $V$ through $y,$
there exists $i\in\{ 0,1,...,n\}$ such that $W=U_\gamma (a_i)$
fulfils $\mu (W)\ge\eta,$ and all chords $yw$ with $w\in W$ form
an angle $\ge\psi$ with $V.$\end{lem}

{\it Proof. } In proving Lemma 4, we found $R,\phi$ such that for
all $y,V$ there is $i$ such that $\dist (a_i,V)\ge R,$ and $ya_i$
and $V$ subtend an angle $\ge\phi .$ Since
$\tan\frac{\phi}{2}<\frac12\tan\phi ,$ we can take
$\gamma=\frac{R}{2}\sin\phi$ in order to ensure that every chord
$yw$ with $w\in U_\gamma(a_i)$ subtends an angle
$\ge\frac{\phi}{2}$ with $V.$ With $\psi=\frac{\phi}{2}$ and
$\eta=\min\{ \mu(U_\gamma(a_i)\,|\, i=0,...,n\},$ the lemma holds
true. \hfill $\Box$\medskip

A {\it self-similar measure} $\mu$ with probability vector
$(p_1,...,p_m)$ with $p_i>0, \sum p_i=1$ and contracting
similarity maps $f_i$ is given by the equation
\[ \mu (B)= \sum_{i=1}^m p_i\mu(f_i^{-1}(B)) \qquad \mbox{ for
Borel sets }B\subset\R^n\ .\] We require $\mu(\R^n)=1.$ Then for
given $p_i,f_i$ there is a unique $\mu ,$ and $\supp\mu$ is the
self-similar set $A$ associated with the $f_i$
\cite{Fal,Bar}.

Let $I=\{ 1,...,m\},$ and let $I^*=\bigcup_{k=1}^\infty I^k$
denote the set of words $\jj=j_1...j_k$ on $I,$ and $I^\infty$ the
set of sequences $j_1j_2...$ We write $f_\jj=f_{j_1}...f_{j_k}$
and $p_\jj=p_{j_1}\cdot ...\cdot p_{j_k}.$ The product measure
$\nu=(p_1,...,p_m)^\infty$ on $I^\infty$ assigns to each cylinder
set $C_\jj=\{ i_1i_2... |\, i_1...i_k=\jj\}$ the value
$\nu(C_\jj)=p_\jj.$ There is the continuous address map $\pi:
I^\infty\to A$ with $\pi(j_1j_2...)=\bigcap_k A_{j_1...j_k}.$  The
measure $\mu$ is the image measure of the product measure,
$\mu=\nu \pi^{-1}$ \cite{Bar}. \vspace{1ex}

{\it Proof of Theorem 1 (ii). } It remains to show the assertion
of the lemma for arbitrary small neighborhoods $U_\epsilon(x)$ by
using self-similarity. We consider the set $\JJ$ of all words
$\jj=j_1...j_k\in I^*$ for which $A_\jj\subset U_\epsilon(x)$ but
$A_{j_1...j_{k-1}}\not\subset U_\epsilon(x).$ Then
\[ \mu (U_\epsilon(x))=\sum_{\jj\in\JJ} p_\jj \] because
$\pi^{-1}(U_\epsilon (x))$ is an open subset of $I^\infty,$ that
is, a countable union of disjoint cylinder sets, which are the
$C_\jj, \jj\in\JJ .$

Now let a hyperplane $V'$ through $x$ be given. For all
$\jj\in\JJ$ we apply $f_\jj^{-1}$ with $y_\jj=f_\jj^{-1}(x)$ as in
the proof of Theorem 1 (i). Lemma 5 yields a set $W_\jj\subset A$
with $\mu(W_\jj)\ge\eta$ such that all chords $y_\jj w$ with $w\in
W_\jj$ subtend an angle $\ge\psi$ with $V_\jj=f_\jj^{-1}(V').$ Now
the set \[ W'=\bigcup_{\jj\in\JJ} f_\jj(W_\jj)\quad \mbox{ fulfils
}\quad \mu(W')\ge\sum_{\jj\in\JJ} \eta p_\jj = \eta\mu
(U_\epsilon(x))\, .\] Note that the measures of $f_\jj(W_\jj)$
have to be added even if the $W_\jj$ overlap. Moreover, all chords
$xw'$ with $w'\in W'$ subtend an angle $\ge\psi$ with $V'.$ This
holds for fixed $x,V'$ with arbitrary $\epsilon .$ Taking
$\delta<\eta ,$ the condition for approximate tangent cannot be
fulfilled. \hfill $\Box$\medskip

\begin{rem} It is not necessary that the mappings $f_\jj$ between
$A$ and the small pieces are similitudes. It is enough to require
that there exists a constant $C$ such that for all $\jj\in I^*$
and all angles $\beta=\angle abc$ between points in $A,$ the image
angle fulfils $f_\jj(\beta)= \angle f_\jj (a)f_\jj(b)f_\jj(c) \ge
C\beta .$ Thus Theorem 1 immediately extends to self-conformal
sets as studied in \cite{K0,K1}.
\end{rem}

\section{Self-affine curves with two pieces}
We shall consider two contracting affine mappings
$f_i(x)=M_ix+v_i$ which are not similarity maps. So the
eigenvalues $\nu_i,\lambda_i$ of $f_i$ must be real, and we
further assume that they are positive:
\qquad $ 0<\nu_i\le\lambda_i<1$\quad for $i=1,2.$\\
The fixed points of the $f_i$ will now be taken as unit points
$e_1={1\choose 0}$ and $e_2={0\choose 1}$ of our coordinate
system. The eigendirections of $M_i$ with respect to $\lambda_i$
are taken as axes. We assume that the eigendirections of
$\lambda_1$ and $\lambda_2$ are not parallel, so that the axes do
intersect.

Compared to the coordinate system of Figure 1, we turned our axes
by 45 degrees. With these coordinates, our mappings have the form
\begin{equation} f_1(x)={\lambda_1\ \alpha\choose 0\,\
\nu_1}\cdot x +{1-\lambda_1\choose 0}\ , \qquad f_2(x)={\nu_2\
0\choose \beta\ \lambda_2}\cdot x +{0\choose 1-\lambda_2}\
.\end{equation} The condition $f_1 {0\choose 1}=f_2{1\choose 0}$
implies $\alpha=\nu_2+\lambda_1-1$ and $\beta=\nu_1+\lambda_2-1,$
so that the vector of eigenvalues
$(\lambda_1,\nu_1,\lambda_2,\nu_2)$ parametrizes all our possible
self-affine curves. Moreover, the coordinates of the point
$z=f_1(e_2)=f_2(e_1)$ are ${\nu_2\choose \nu_1}.$

If $\nu_1+\nu_2=1,$ then $J$ will be the segment with endpoints
$e_1,e_2.$ In the case $\nu_1+\nu_2>1$ the set $J$ will not be a
simple curve. Among others, $z$ will be a multiple point of the
curve. Thus we shall require $\nu_1+\nu_2<1,$ which means that $z$
lies in the interior of the triangle $T$ with vertices $0,e_1,$
and $e_2.$ Note that $\lambda_1+\nu_2<1$ in Theorem 2 implies
$\nu_1+\nu_2<1.$ Let us also note that in the special case
$\nu_1=\lambda_1$ we get $\alpha<0.$ Similarly, $\nu_2=\lambda_2$
implies $\beta<0$ so that we have no similarity maps.
\medskip

\noindent {\bf Assumptions. }\emph{The conditions of Theorem 2 are
given in (1) above and (2) below. They are taken as assumptions
for the rest of the paper.}
\begin{equation}
0<\nu_i\le \lambda_i<1 \ \mbox{ for }i=1,2,\quad
\nu_1+\nu_2<1,\quad \alpha=\nu_2+\lambda_1-1, \quad
\beta=\nu_1+\lambda_2-1\smallskip
\end{equation}

\begin{figure}\hspace*{-2.5cm}
\includegraphics[width=.85\textwidth]{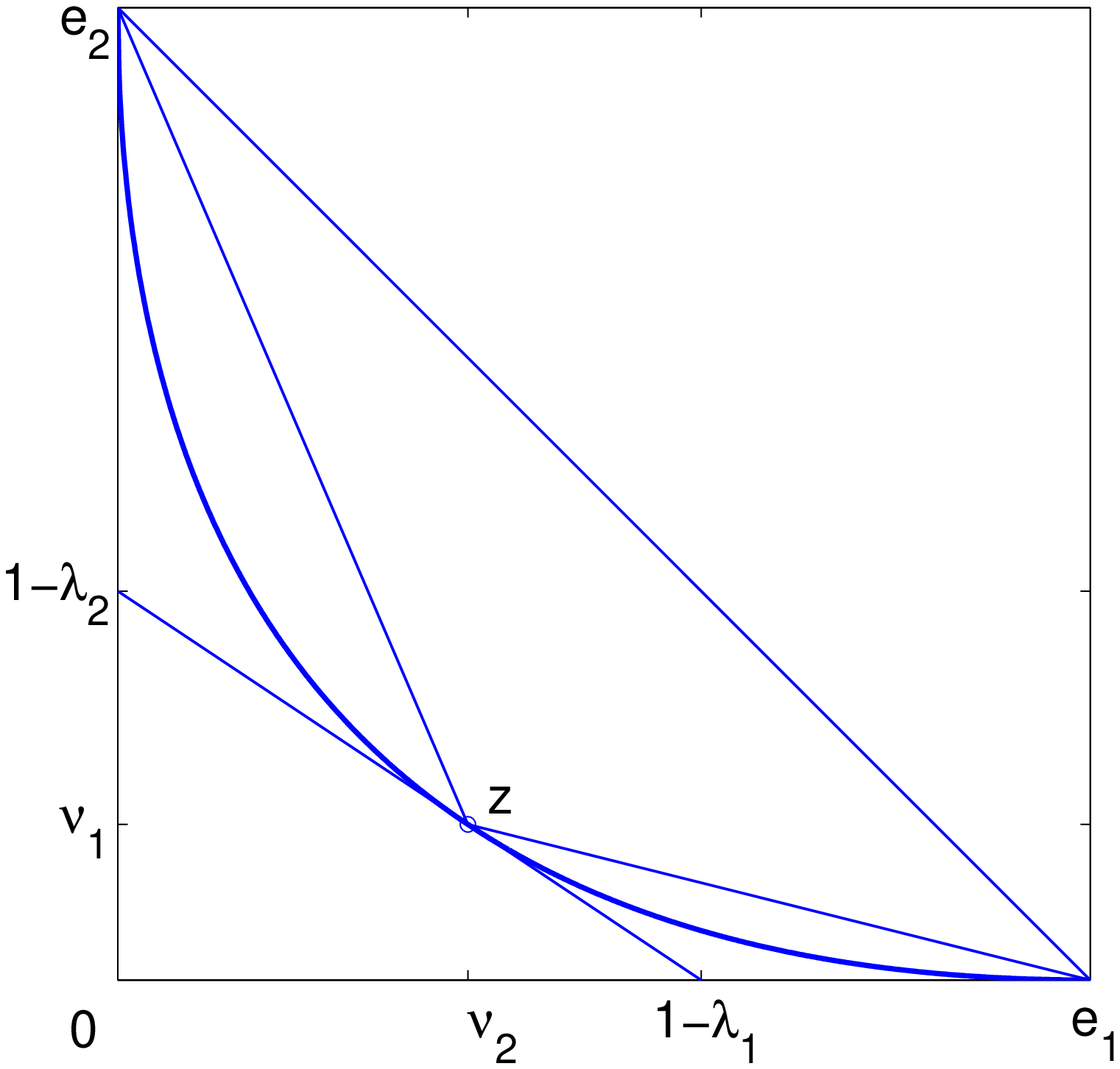}\hspace*{-2cm}
\includegraphics[width=.41\textwidth]{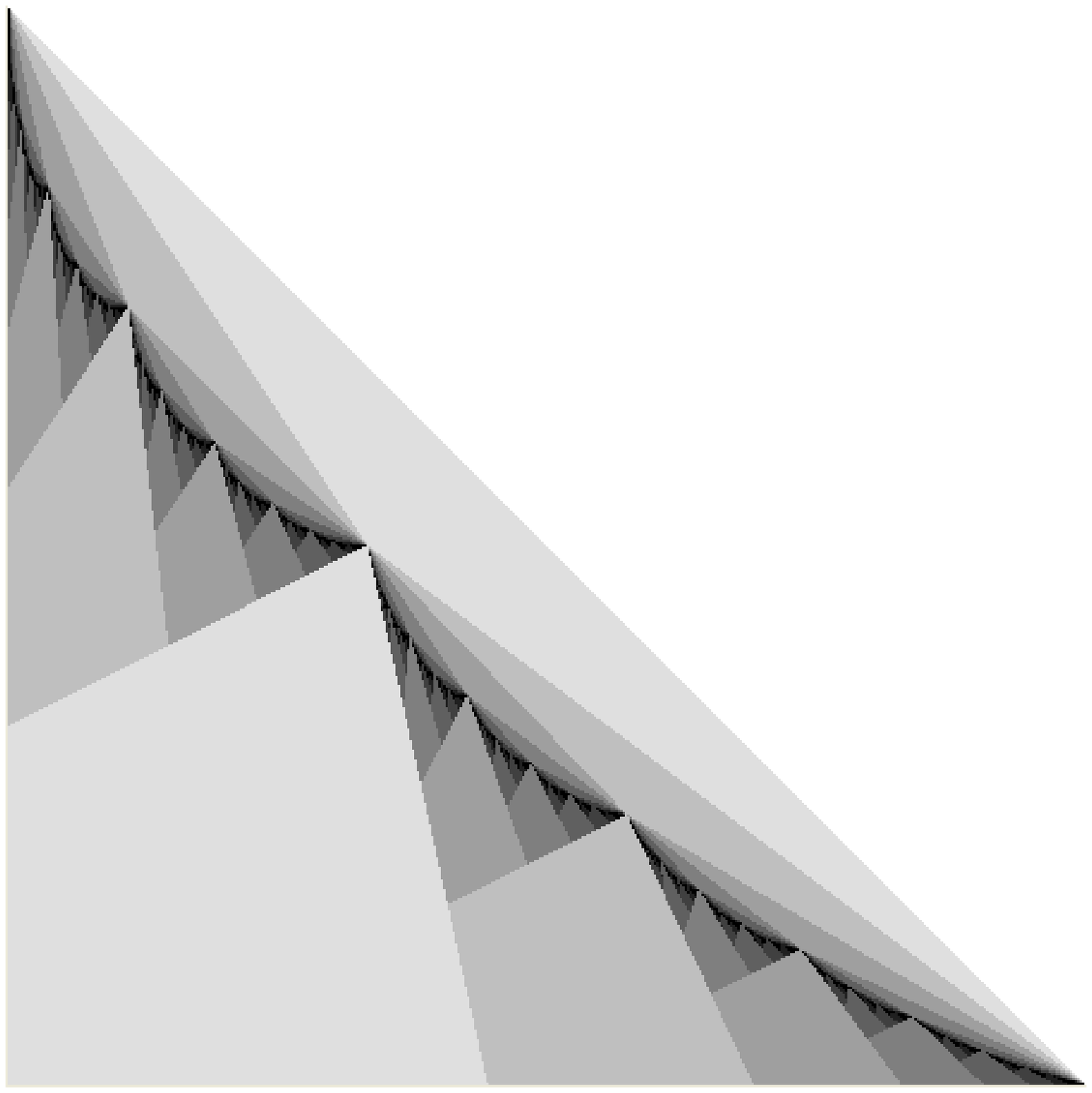}
\vspace{1.5ex}\\
\caption{A smooth curve ($\alpha, \beta <0$) and a
non-differentiable curve ($\alpha<0,\beta>0$) in the new
coordinate system. The triangle $T=\triangle 0e_1e_2$ and its
images under the $f_i$ are indicated.}
\end{figure}

The maps $f_1,f_2$ are said to fulfil the \emph{open set
condition} if there is an open set $U$ such that $f_1(U)$ and
$f_2(U)$ are disjoint subsets of $U$ \cite{Mo,Fal,Bar}. For
self-similar curves, this condition is hard to verify even if we
assume that $f_1(J)\cap f_2(J)$ consists of a single point only
\cite{ATK,BR}. In our setting, we get this condition for free.

\begin{Prop} The intersection of $f_1(J)$ and $f_2(J)$ consists of a single point. If $U$
denotes the interior of the triangle $T$ with vertices ${0\choose
0}, {1\choose 0}, {0\choose 1}$ then $f_1,f_2$ fulfil the open set
condition with $U.$\end{Prop}

\begin{Proof} We just note that $f_1(T)$ is the triangle with vertices
$e_1,z,{1-\lambda_1\choose 0},$ and $f_2(T)$ is the triangle with
vertices $e_2,z,{0\choose 1-\lambda_2}.$ See Figure 2. \end{Proof}

\begin{Prop} $J$ is a simple curve without double points. If for the moment we
write $f_0$ instead of $f_2,$ an explicit homeomorphism $\varphi:
[0,1]\to J$ is given by binary numbers: $\varphi(t)=
\bigcap_{k=1}^\infty f_{i_1...i_k}(J)$ for $t=\sum_{k=1}^\infty
i_k\cdot 2^{-k},\ i_1,i_2,...\in\{ 0,1\} .$ \end{Prop}

\begin{Proof} This directly follows from Proposition 7. See
\cite{Hut}, 3.5, compare also \cite{ATK,BR}, or \cite{Bar},
Chapter VIII.2. The point is that the two addresses
$01111...=0\overline{1}$ and $1\overline{0}$ of the intersection
point $z=J_1\cap J_2$ are the same as for the point $\frac12$ when
$[0,1]$ is considered as self-similar set with respect to
$g_0(t)=\frac{t}{2}$ and $g_1(t)=\frac{t+1}{2}.$\end{Proof}

\begin{Prop} If \qquad $\lambda_1+\nu_2\le 1$\quad and \quad $\lambda_2+\nu_1\le
1$\\ then $J$ is the graph of a strictly decreasing function
$x_2=\psi (x_1)$ in our coordinate system. Thus $J$ is
rectifiable, with one-dimensional Hausdorff measure $0<\HH^1(J)\le
2.$ The function $\psi$ is differentiable at almost all points
$x_1.$
\end{Prop}
\begin{Proof}
Let  $W=\{ {w_1 \choose w_2}|\, w_1w_2\le 0\}$ denote the cone
consisting of the second and fourth quadrant of our coordinate
system. Our assumptions say that $\alpha, \beta\le 0,$ so $W$ is
invariant under the linear mappings $h_i(x)=M_ix$ corresponding to
the $f_i$ in (1): $M_iw\in W$ for $w\in W.$

The curve $J$ is the graph of a decreasing function if $y-y'\in W$
for all $y,y'\in J.$ (Then $y-y'$ cannot be a multiple of $e_2$
(or $e_1$, respectively) because in such case all points between
$y$ and $y'$ must lie on that vertical segment which contradicts
self-affinity. See also Section 4.) Since $W$ is closed, and any
$y\in J$ can be approximated by points of the form $f_\jj(e_i)$
with $i\in I$ and $\jj\in I^*,$ it suffices to show the relation
only for points $y,y'$ of this form.

Clearly $e_2-e_1$ is in $W.$ For $x-x'\in W$ and $i=1,2$ we have
$f_i(x)-f_i(x')=M_i(x-x')\in W.$ This implies
$v=f_\jj(e_2)-f_\jj(e_1)\in W$ for every $\jj\in I^*.$ For
arbitrary $k\in\N$ and $\ii,\jj\in I^k,$ the difference
$f_\jj(e_2)-f_\ii(e_1)$ is a finite sum of vectors of the form
$v,$ since in the $k$-th level of construction of $J,$ the
triangles $f_\jj(T)$ and $f_\ii(T)$ are connected by a finite
chain of such triangles. Since finite sums of vectors in $W$
belong to $W,$ we proved $f_\jj(e_2)-f_\ii(e_1)\in W$ for
$\ii,\jj\in I^k$ and $y-y'\in W$ for $y,y'\in J.$

Since the curve $J$ is the graph of a function, its length is
\begin{equation} \HH^1(J)=\sup\{\sum_{i=1}^N |y^{(i)}-y^{(i-1)}|
\, |\, N\in\N , 0=y_1^{(0)}<y_1^{(1)}<...<y_1^{(N)}=1 \}\
.\end{equation} Since $|y^{(i)}-y^{(i-1)}|\le
|y_1^{(i)}-y_1^{(i-1)}|+|y_2^{(i)}-y_2^{(i-1)}|$ and $\sum_{l=1}^N
|y_1^{(i)}-y_1^{(i-1)}|=\sum_{l=1}^N |y_2^{(i)}-y_2^{(i-1)}|=1,$
we get $\HH^1(J)\le 2.$ The inequality $\HH^1(J)>0$ holds since
$J$ is connected, see \cite[Lemma 3.2]{Fa}. It is well-known that
curves of finite length are differentiable $\HH^1$-almost
everywhere \cite[Theorem 3.8]{Fa}, and monotone functions are
differentiable at almost every real argument.
\end{Proof}

\begin{figure}
\includegraphics[width=.46\textwidth]{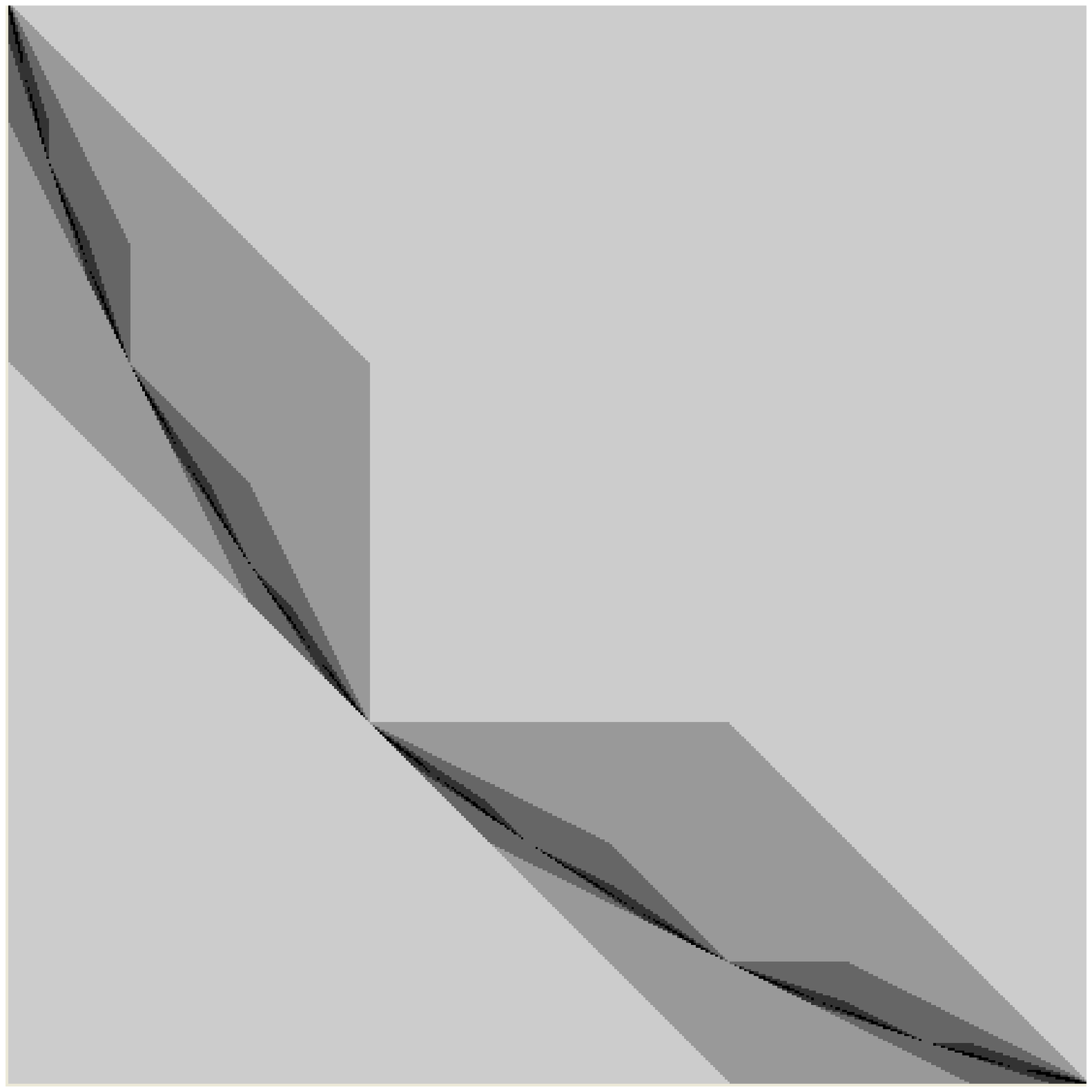}
\includegraphics[width=.46\textwidth]{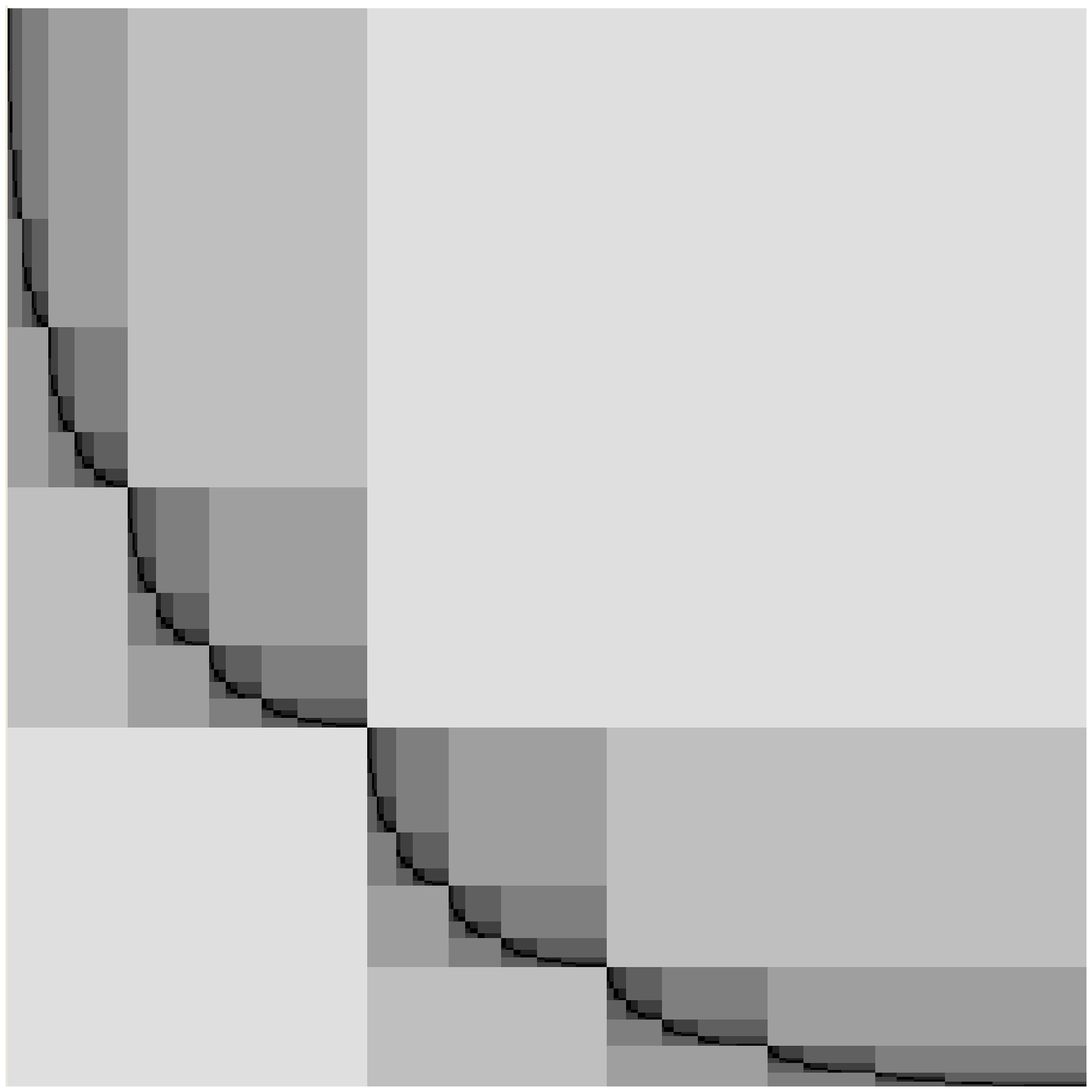}
\vspace{1.5ex}\\
\caption{The images of the unit square under the $f_\ii.$ The
curve on the left is smooth, $\lambda_i=\nu_i=\frac{1}{3}.$
Example 10 on the right has uncountable dense sets of points where
the curve is differentiable and non-differentiable,
$\lambda_i=\frac{2}{3}, \nu_i=\frac{1}{3}.$ }
\end{figure}

\begin{Exa} The case $\lambda_1=\lambda_2=\frac{2}{3},
\nu_1=\nu_2=\frac{1}{3}, \alpha=\beta=0,$ illustrated in Figure
1b, was studied by K\"aenm\"aki and Vilppolainen \cite[Example
6.2]{KV}. $J$ fulfils the conditions of Proposition 9 and thus has
a tangent $\HH^1$-almost everywhere although at all points
$y=f_\jj(z)$ the left-sided derivative is 0 and the right-sided
one is $\infty .$
We have $\HH^1(J)=2 .$\\
Let $\mu$ be the self-similar measure on $J$ induced by the
product measure $(\frac12,\frac12 )^\infty$ (see Section 2). Then
at $\mu$-almost all points, $J$ does \emph{not} possess a tangent.
\end{Exa}

\begin{Proof} We consider the rectangles $R_\jj=f_\jj([0,1]^2)$
for $\jj\in I^n,$ with side length $s_\jj$ in direction $e_1$ and
$t_\jj$ in direction $e_2.$ The ratio $r_\jj=s_\jj/t_\jj$ of the
side lengths is $2^{n_\jj}$ where $n_\jj=\# \{ k\le n\,|\, j_k=1\}
-\# \{ k\le n\,|\, j_k=2\}.$ Thus $-n\le n_\jj\le n .$
Probabilistic arguments will now show that for large $n,$ most of
the rectangles are very lengthy and slim.

When all words $\jj\in I^n$ are assigned equal probability
$2^{-n},$ then $n_\jj$ is distributed like the endpoints of a
symmetric random walk on the integers, with start in zero, after
$n$ steps.

The two rectangles at $f_{i_1...i_k}(z)$ correspond to the words
$\jj=i_1...i_k122...2$ and $\jj '=i_1...i_k211...1.$ Thus
$n_\jj\to-\infty$ and $n_{\jj '}\to\infty$ for $n\to\infty,$ which
shows the assertion on one-sided derivatives.\smallskip

Now we shall show that for large $n,$ the sum of $\diam R_\jj$
over $\jj\in I^n$ is almost 2, which by (3) implies $\HH^1(J)\ge
2.$ We let $n$ be even. Then ${\bf P}(n_\jj=0)\approx
\frac{1}{\sqrt{\pi n}}$ is the middle term of the binomial
distribution for $n$ and $p=\frac12,$  so ${\bf P}(n_\jj\le
0)<\frac12+\frac{1}{\sqrt{n}}$ holds for large $n$ (in fact for
all even $n$).

For $n_\jj\le 0$ the square has the biggest horizontal side:
$s_\jj<(\frac{2}{9})^{n/2}.$ Thus \[ \sum_{n_\jj\le 0} s_\jj \le
2^n(\frac12+\frac{1}{\sqrt{n}})(\frac{2}{9})^{n/2}\le
(\frac{8}{9})^{n/2}\to 0 \quad\mbox{ for } n\to\infty\ .\] This
implies that $\sum_{n_\jj\ge 1} \diam R_\jj\ge \sum_{n_\jj\ge 1}
s_\jj\to 1.$ A similar calculation for the $t_\jj$ yields
$\sum_{n_\jj\le -1} \diam R_\jj\to 1$ for $n\to\infty .$ So we
obtain $\HH^1(J)\ge 2.$\smallskip

For the last statement, we use the recurrence of the
one-dimensional symmetric random walk. For almost all
$j_1j_2...\in I^\infty,$ with respect to the product measure
$(\frac12,\frac12 )^\infty,$ we have $n_{j_1...j_n}=0$ for
infinitely many $n.$ That is, $\mu$-almost all points $y\in J$ lie
in arbitrarily small squares $R_\jj.$ These points $y$ have
arbitrary small neighborhoods which are geometrically similar to
$J.$ Using the argument of the proof of Theorem 1 in Section 2, we
can conclude that $J$ does not admit a tangent at such points $y.$
\end{Proof}

In this example, the points of differentiability form an
uncountable dense set with an uncountable dense complement. It
turns out, however, that all tangent lines have only two
directions, horizontal and vertical (cf. Section 4). In \cite{KV}
the example was used to show that there are exactly two
equilibrium measures. These measures are now obtained as image
measures of the Lebesgue measure on the axes under the natural
one-to-one map. They assign $s_\jj$ and $t_\jj,$ respectively, to
$R_\jj.$ The Hausdorff measure is just the sum of these two
measures.

\section{Smooth self-affine curves}
\begin{Prop} The curve $J$ has the two coordinate axes
as tangents at $e_1$ and $e_2,$ respectively. Thus $J$ has
one-sided tangents at all points where two little pieces meet.
\end{Prop}

\begin{Proof} As neighborhoods of $e_1$ in $J$ we consider the
sets $f_1^k(J)$ for large $k.$ Since $f_1^k(J)\subset f_1^k(T),$
we just have to consider the angle $\beta_k$ of the triangle
$f_1^k(T)$ at the vertex $e_1$ (cf. Section 1). From (1) we have
\[ f_1^k(x)={\lambda_1^k\ \alpha_k\choose 0\,\
\nu_1^k}\cdot x +{1-\lambda_1^k\choose 0}\ . \] In the case
$\lambda_1=\nu_1,$ induction shows
$\alpha_k=k\alpha\lambda_1^{k-1}$ where
$\alpha_1=\alpha=\nu_2+\nu_1-1<0.$ For $\lambda_1>\nu_1,$
induction shows
\[\alpha_k=\alpha\cdot\frac{\lambda_1^k-\nu_1^k}{\lambda_1-\nu_1}\
.\] The relation $\alpha=\nu_2+\lambda_1-1<\lambda_1-\nu_1$ can be
rewritten as $\alpha=\gamma(\lambda_1-\nu_1)$ for some constant
$\gamma<1.$ This implies
$\alpha_k=\gamma\cdot(\lambda_1^k-\nu_1^k).$

The vertices of the triangle $f_1^k(T)$ are $e_1$ and
$f_1^k(0)=(1-\lambda_1^k)e_1$ and
\[ f_1^k(e_2)={\alpha_k+1-\lambda_1^k\choose \nu_1^k}={y_1\choose y_2}\ .\]
For the angle $\beta_k$ at $e_1$ we have $\tan
\beta_k=y_2/(1-y_1)=\nu_1^k/(\lambda_1^k-\alpha_k).$ If
$\lambda_1=\nu_1$ then $\tan \beta_k=\lambda_1/(\lambda_1-\alpha
k)$ with $\alpha<0.$ For the case $\lambda_1>\nu_1$ we get $\tan
\beta_k=\nu_1^k/(\lambda_1^k(1-\gamma)+\nu_1^k\gamma).$ In both
cases, $\beta_k$ converges to 0.

For $e_2,$ the proof is similar. Now consider
$z=f_1(e_2)=f_2(e_1).$ Since affine maps preserve tangents of
curves, $f_1([0,e_2])=[z, (1-\lambda_2)e_2]$ is a right-sided
tangent of $J$ at $z,$ and $f_2([0,e_1])=[z, (1-\lambda_1)e_1]$ is
a left-sided tangent at $z.$ By the same argument, all points
$f_\jj(z)$ with $\jj\in I^*$ admit one-sided tangents. These are
the points where two pieces of $J$ meet.
\end{Proof}

\noindent {\it Proof of Theorem 2.\  (i) } We show that if
$\alpha=\lambda_1+\nu_2-1<0$ and $\beta=\lambda_2+\nu_1-1<0,$ then
$J$ admits tangents at all points which are not endpoints of
little pieces.  This is done with Proposition 12 below.

{\it (ii) } Since $z={\nu_2 \choose \nu_1},$ Proposition 11 says
that the one-sided tangents at $z$ have slope
$(1-\lambda_2-\nu_1)/\nu_2$ and $\nu_1/(1-\lambda_1-\nu_2),$
respectively. If $\nu_1\nu_2=
(1-\lambda_1-\nu_2)(1-\lambda_2-\nu_1)$ then the one-sided
tangents agree, and $J$ admits tangents at $z$ and all $f_\jj(z).$
Thus $J$ is a differentiable curve. The continuity of the
derivative is shown in Remark 15 below.  {\ \hfill $\Box$\medskip}

\begin{Prop} If $\alpha<0$ and $\beta<0$ then the curve $J$ has
tangents at all points which are not endpoints of pieces.
\end{Prop}

\begin{Proof} For Proposition 9, we proved that all chords $y-y'$
with $y,y'\in J$ lie in the cone $W=\{ {w_1 \choose w_2}|\,
w_1w_2\le 0\}$ given by the second and fourth quadrant. For the
matrices $M_i$ of $f_i$ we have $M_iW\subseteq W, i=1,2.$  Thus
the chords between points of $J_{i_1...i_n}$ lie in
$M_{i_1}...M_{i_n}W\subseteq W.$ Now we show that for each
sequence $\ii=i_1i_2...\in I^\infty$ the cones
$M_{i_1}...M_{i_n}W$ converge to a single line $L=L(\ii)$ for
$n\to\infty .$ Then $L+x$ is the tangent of $J$ at $x=\pi(\ii).$

For sequences $\ii$ which end with $\overline{1}=111...$ or with
$\overline{2}$ this can be shown as in Proposition 11 (see Remark
14). But we obtain only one-sided tangents since $x$ was an
endpoint of $J_{i_1...i_n}$ for all sufficiently large $n.$ Now we
consider an arbitrary point $x=\pi(\ii)$ which is not an endpoint
of any small piece of $J,$ and prove there is a two-sided
derivative. Since $\ii$ does not end with $\overline{1}$ or
$\overline{2},$ there exist positive integers $1<k_1<k_2<...$ such
that $i_{k_m}=1$ and $i_{k_m+1}=2$ for $m=1,2,...$

The size of a subcone $W'\subseteq W$ is taken as $\diam
([-e_1,e_2]\cap W').$ It will be helpful to define the linear
functional $\phi (w)=w_2-w_1,$ and take the line $\ell =\{ w\,|\,
\phi(w)=1\}$ instead of $[-e_1,e_2].$ We shall find a constant
$\delta>0$ such that
\begin{equation}
\diam (\ell\cap M_{i_{k_1}}...M_{i_{k_m-1}}W) \le
(1-\delta)^{m-1}\diam (\ell\cap W)
\end{equation}
for all $m.$ Then $M_{i_{k_1}}...M_{i_{k_m-1}}W$ converges to a
line for $m\to\infty ,$ and the proof is complete.

Let $V=M_1M_2W.$ Since $\alpha, \beta <0$ and
\[ M_1 M_2=
 \left(\begin{array}{cc} \lambda_1 & \alpha \\ 0 & \nu_1 \end{array}\right)
 \left(\begin{array}{cc} \nu_2 & 0 \\ \beta & \lambda_2 \end{array}\right)=
 \left(\begin{array}{cc} \lambda_1\nu_2+\alpha\beta & \alpha\lambda_2 \\ \nu_1\beta & \nu_1\lambda_2 \end{array}\right).
 \]
both $M_1M_2(-e_1)$ and $M_1M_2e_2$ are in $\inte W$ and hence
$V\subset\inte W.$ Consequently,
\[ \delta:=\frac{\dist(\ell\cap V, \ell\cap\partial W)}{\diam(\ell\cap W)}
=\frac{\dist([-e_1,e_2]\cap V, \{ -e_1,e_2\})}{\diam [-e_1,e_2]}>0
\ ,
\]
and $\diam (\ell\cap V)\le (1-2\delta)\cdot\diam (\ell\cap W)\, .$

\begin{lem} If $B$ is a regular $2\times 2$ matrix with
$BW\subseteq W$ then
\[ \diam (\ell\cap BV)\le (1-\delta)\cdot\diam (\ell\cap BW)\, .
\]
\end{lem}

\noindent {\it Proof of Lemma. \ } Let $\ell'=B(\ell).$ Since
linear maps preserve length ratios within lines, we have
\begin{equation} \delta=\frac{\dist(\ell'\cap BV, \ell'\cap B(\partial
W))}{\diam(\ell'\cap BW)}
\end{equation}
$\ell'\cap B(\partial W)$ consists of two points $-Be_1, Be_2$
which we call $a$ and $b,$ choosing the order so that
$\phi(a)/\phi(b)\ge 1$. Furthermore, let $c$ be the point in
$\ell'\cap BV$ which assumes maximum distance from $a.$ Then
$c=\tau a+(1-\tau)b$ with $\tau\in[\delta, 1-\delta]$ because of
(5).

The projection $p(w)=w/\phi(w)$ maps $[a,b]$ onto $\ell\cap BW,$
and $\ell\cap BV$ is contained in $[p(a),p(c)].$ Thus
\begin{eqnarray*} \diam(\ell\cap BV)&\le &|p(c)-p(a)|=
  \left|\frac{\tau a+(1-\tau)b}{\tau \phi(a)+(1-\tau)\phi(b)}-\frac{a}{\phi(a)}\right|
  \\
 &=&\left|\frac{\tau a+(1-\tau)b - \tau a -
(1-\tau)a\cdot\phi(b)/\phi(a)}{\tau
\phi(a)+(1-\tau)\phi(b)}\right|\\
 &=&\frac{(1-\tau)\cdot|b/\phi(b)-a/\phi(a)|}{\tau \phi(a)/\phi(b)+(1-\tau)}
 \\ &\le & (1-\delta)\cdot|p(a)-p(b)| =(1-\delta)\cdot\diam(\ell\cap BW)
 \ . \qquad\qquad\quad \Box
\end{eqnarray*}

To finish the proof of Proposition 12, we show (4) by induction on
$m.$ Let $B_m=\prod_{k=1}^{k_m -1} M_{i_k}.$ The matrix
$C_m=\prod_{k=k_m}^{k_{m+1}-1} M_{i_k}$ has the form
$C_m=M_1M_2\cdots$, so $C_m W\subset M_1M_2W=V$ and $B_{m+1} W=B_m
C_m W\subset B_m V$ for $m=1,2,...$ \  We apply the Lemma to
$B=B_m$ and obtain
\[ \diam(\ell\cap B_{m+1} W)\le\diam(\ell\cap B_m V)
\le(1-\delta)\cdot\diam(\ell\cap B_m W) \]
for $m\ge 1.$ This implies (4).
\end{Proof}

\begin{rem} Assuming $\alpha,\beta <0,$ we proved that for each address sequence $i_1i_2...$
which does not end with $\overline{1}$ or $\overline{2},$ and the
corresponding point $x=\pi(i_1i_2...)$ there is the tangent line
$L(x)=\bigcap_{n=1}^\infty M_{i_1...i_n}W\, .$ This is also true
for $i_1i_2...=\overline{1}$ or $\overline{2},$ hence for all
addresses.
\end{rem}

\begin{Proof} The cones $M_1^n(W)$ form a decreasing sequence, and
their intersection is a closed convex cone in $\R^2.$ If this is
not the line $r e_1, r\in \R,$ the intersection cone has a second
bounding line determining another eigenvector of $M_1.$ However,
the second eigenvector ${-\alpha\choose \lambda-\nu}$ of $M_1$
does not belong to $W.$ Note that $\alpha<0$ is needed here. In
Proposition 11 we used a subcone of $W.$ The proof for
$\overline{2}$ uses $\beta<0,$ and
$M_{i_1...i_n\overline{1}}=M_{i_1...i_n}L(e_1)$ is a line.
\end{Proof}

\begin{rem}
For each address $i_1i_2...$ of $x\in J$ the $M_{i_1...i_n}W$
shrink down to a line $L(x).$ Only when $i_1i_2...$ ends with
$\overline{1}$ or $\overline{2},$ there is a second address of $x$
and a second line $L'(x).$ So for each $\epsilon>0$ there is $n$
such that $\diam(\ell\cap M_{i_1...i_n}W)<\epsilon .$ Thus all
tangents $L(y)$ (and also $L'(y)$) with $y\in J_{i_1...i_n}$ are
within distance $\epsilon$ from $L(x).$ If $x$ is not an endpoint
of a piece of $J,$ the $J_{i_1...i_n}$ are neighborhoods of $x,$
and the derivative is continuous at $x.$

For endpoints, we have one-sided neighborhoods, and one-sided
continuity. If now the one-sided tangents at $z$ coincide, then
$L(x)=L'(x)$ for all points with two addresses, we can put the
one-sided neighborhoods together, and the derivative becomes a
continuous function.
\end{rem}

\section{The uniqueness of the parabola}
One can ask whether our $C^1$ self-affine functions are even
$C^2.$ There is at least one well-known example, the parabola.
Parabolic arcs play the same role for self-affine sets as
intervals do for self-similar sets. The parabola $x_2=x_1^2/2$
admits a transitive group of `parabolic translations' $\phi({x_1
\choose x_2})={x_1+t\choose tx_1+x_2+t^2/2}, \ t\in \R ,$ as well
as `parabolic homotheties' $\psi({x_1 \choose x_2})={\lambda
x_1\choose \lambda^2 x_2}$ with $\lambda>0.$ Thus any arc $P$ of a
parabola can be mapped by an affine map into any other one, and
for every division $P=P_1\cup P_2$ into two subarcs, generated by
an interior point $z$ of $P,$ we have a representation of $P$ as a
self-affine set.

\begin{Exa} In our coordinate system, the symmetric representation
of a parabolic arc $P$  is given by $\lambda_i=\frac12,
\nu_i=\frac{1}{4}.$ The following one-parameter family provides
all representations of $P $ (cf.~Figure 2 for $\lambda_1=0.4$).
\[ \lambda_1+\lambda_2=1\ ,\quad  \nu_1=\lambda_1^2\ ,\quad \nu_2=\lambda_2^2 \
\]
\end{Exa}

\begin{Proof}
A parabola is defined by two of their points together with their
tangents, so $P$ is uniquely determined. In the coordinate system
of Figure 1 it has equation $x_2=\frac12(x_1^2+1).$ This is
$2w_1+2w_2-1=(w_1-w_2)^2$ in our coordinate system since
$x_1=w_1-w_2, x_2=w_1+w_2 .$ Solving the quadratic equation for
$w_2$ we get
\[ w_2=w_1+1-2\sqrt{w_1}=(1-\sqrt{w_1})^2\quad , \qquad\mbox{ or
}\qquad \sqrt{w_1}+\sqrt{w_2}=1\ .
\]
The derivative is $w_2'=1-1/\sqrt{w_1}.$ Now consider the tangent
line to $P$ at an arbitrary point ${\nu_2\choose \nu_1}.$ It has
equation $w_2=\nu_1+ (w_1-\nu_2)(1-1/\sqrt{\nu_2}).$ This line
passes through $w_1=0, w_2=1-\lambda_2$ which yields
\[ 1-\lambda_2=\nu_1-\nu_2+\sqrt{\nu_2}=1-\sqrt{\nu_2}\quad\mbox{
because of }\quad \nu_1=1+\nu_2-2\sqrt{\nu_2} \ .
\]
Thus $\lambda_2=\sqrt{\nu_2},$ and similarly
$\lambda_1=\sqrt{\nu_1},$ since the tangent line also passes
through $w_1=1-\lambda_1, w_2=0.$ And
$\sqrt{\nu_1}+\sqrt{\nu_2}=1$ was already proved.
\end{Proof}

\noindent {\it Proof of Theorem 3, (i). \ } We have a plane $C^2$
curve $J$ parametrized by $x=\phi(t)$ with $\phi(0)=0$ and
$\phi''(0)\not= 0,$ and a contractive linear map $h\not= 0$ which
maps $J$ into itself. We show that $J$ is contained in a parabola.
The condition $\phi''(0)\not= 0$ does not depend on $\phi .$ It
says that the curvature of $J$ at 0 is nonzero. In particular, any
subarc of $J$ with the point 0 cannot be contained in a line.

Since $h(J)\subseteq J$ and $h(0)=0,$ the map $h$ fixes the
tangent line of $J$ at 0. This must be an eigenvector of $h.$ We
use the Jordan normal form of $h.$ The classes of $C^1$ and $C^2$
curves, as well as the class of parabolas, and a curve's property
of having nonzero curvature at a point, are preserved under linear
transformations. For that reason, we can \emph{calculate with the
Jordan base as standard coordinate system.} Moreover, we assume
that the first basis vector is an eigenvector of $h$ in the
direction of the tangent of $J$ at 0. Thus $h$ has a matrix of the
form ${\lambda\ \gamma\choose 0\ \nu }$ where $\gamma=0$ if
$\lambda\not= \nu .$ Note that $|\lambda|,|\nu|<1$ since $h$ is
contractive.

In the following, we use coordinates $x,y$ instead of $x_1,x_2.$
Let us take a $C^2$ parametrization $\varphi (t)={x(t)\choose
y(t)}$ of $J$ in our new coordinate system with $\varphi (0)=0$
and $\varphi' (0)\not= 0$ (the parametrization by arc length
fulfils this condition). Then $x'(0)\not= 0$ since $y'(0)=0$ by
the choice of the first axis. So there exists a $\delta>0$ such
that $x'(t)\not= 0$ for $0\le t<\delta .$ Thus $x(t)$ is a
monotone $C^2$ function for $0\le t<\delta .$

By the inverse function theorem, the inverse function $t(x)$ is
$C^2$ for $0\le x<\epsilon =x(\delta) ,$ and the composition
$y(t(x))$ is also $C^2$ for $0\le x<\epsilon .$ This function
which describes $J$ in the vicinity of 0 will now be called
$y=y(x), $ and $t$ is eliminated. We have $y(0)=0,\ y'(0)=0,$ and
$y''(0)\not= 0, $ since otherwise the curvature of $J$ at 0 would
be zero. Now we apply the assumption $h(J)\subseteq J.$ From
\[ {\lambda\  \gamma\choose 0\  \nu }\cdot {x\choose y(x)}=
 {\lambda x+\gamma y(x)\choose \nu y(x)}={z\choose y(z)}
\]
with $z=z(x)=\lambda x+\gamma y(x)$ we get the derivatives
\[ \nu y'(x)= y'(z)\cdot z'(x)\quad \mbox{ with }z'(x)=\lambda+\gamma y'(x)\quad \mbox{ and }
\]
\begin{equation} \nu y''(x)= y''(z)\cdot z'(x)^2+y'(x)\cdot z''(x)\quad \mbox{ with }
z''(x)=\gamma y''(x)
\end{equation}
for $0\le x<\epsilon .$ Inserting $x=0$ and $z=0, y'=0,\
z'=\lambda $ in (6) we obtain $\nu=\lambda^2$ since $y''(0)\not=
0.$ Thus $\gamma=0$ since $\nu$ and $\lambda$ are different.
Equation (6) now simplifies:
\[  y''(x)=  y''(\lambda x)\quad \mbox{ and hence }\  y''(x)=y''(\lambda^n x)
\]
for $x\in [0,\epsilon )$ and $n=1,2,...$ For the limit $n\to
\infty$ we obtain $y''(x)=y''(0)$ since $|\lambda| <1$ and $y''$
is continuous. Thus $y(x)$ is a quadratic function for $0\le
x<\epsilon .$ Applying $h^{-1}$ several times to that part of $J,$
we see that the arc $J$ is contained in a parabola. {\ \hfill
$\Box$\medskip}

\noindent {\it Proof of Theorem 3, (ii). \ } If $J=f_1(J)\cup
...\cup f_m(J)$ is a self-affine curve and $J_\ii$ with
$\ii=i_1...i_n$ is a small piece of $J,$ then $f_\ii$ is a
contracting affine map with fixed point in $J_\ii$ which maps $J$
into itself.

Now assume $J\subset R^2$ is a $C^2$ curve and not subset of a
line, so that $\phi'' (x)\not= 0$ for some parametrization of $J$
and some $x\in J.$ Since $\phi''$ is continuous, this implies
$\phi'' (y)\not= 0$ for all $y$ in a small piece $J_\ii$
containing $x.$ Shifting our coordinate system to the fixed point
of $f_\ii=h$ and applying part (i), we see that $J$ is contained
in a segment or parabola. {\ \hfill $\Box$\medskip}

\end{document}